\documentclass[AMA,STIX1COL]{./WileyNJD-v2}
\usepackage{moreverb}

\newcommand\BibTeX{{\rmfamily B\kern-.05em \textsc{i\kern-.025em b}\kern-.08em
T\kern-.1667em\lower.7ex\hbox{E}\kern-.125emX}}

\articletype{Research Article}%

\received{\today}
\revised{<day> <Month>, <year>}
\accepted{<day> <Month>, <year>}

\usepackage[scale=1.0]{OldStandard}
\usepackage{libertinust1math}

\usepackage{./mylatexcommands}
\usepackage{./submath}

\newcommand{\rosso}{\color{black}} 
\newcommand{\blue}{\color{black}}   
\newcommand{\erase}[1]{\!}

\newcommand{\Fun} [1]{F(X_{#1},\bv{A})}
\newcommand{\fun}[1]{f(x_{#1},\bv{a})}
\newcommand{\funuw}[1]{f(x_{#1},\bv{a}_{uw})}
\newcommand{\funb}[1]{f(x_{#1},\bv{a})}
\newcommand{\Lpsev}{P_{7}(x)}

\newcommand{\mse}{e^{\ssup{0.6}{2}}}
\newcommand{\msem}{\overline{\mse}}

\newcommand{\Xmin}{X_{min}}
\newcommand{\Xmax}{X_{max}}
\newcommand{\Ymin}{Y_{min}}
\newcommand{\Ymax}{Y_{max}}

\newcommand{\figdir}{.}         

\renewcommand {\Rfi} [1] {\mbox{Figure \ref{#1}\hspace*{-0.5em}}}   
\renewcommand {\Rta} [1] {\mbox{Table \ref{#1}\hspace*{-0.5em}}}    

\newboolean{refcomm}

\setboolean{refcomm}{false}

\ifthenelse{\boolean{refcomm}}
{%
 \newcommand{\refereeone}{\color{blue}}%
 \newcommand{\refereetwo}{\color{magenta}}%
}
{
 \newcommand{\refereeone}{\color{black}}%
 \newcommand{\refereetwo}{\color{black}}%
}

\begin{document}

\title{Maximal-entropy driven determination of weights in least-square approximation}

\author[1a]{Domenico Giordano*}

\author[2]{Felice Iavernaro}


\authormark{Domenico Giordano \textsc{et al}}

\address[1]{\orgdiv{European Space Agency}, \orgname{ESTEC (retired)}, \orgaddress{\country{The Netherlands}}}

\address[2]{\orgdiv{Dipartimento di Matematica}, \orgname{Universit\`a di Bari}, \orgaddress{\country{Italy}}}


\corres{*\email{dg.esa.retired@gmail.com}}

\presentaddress{c/o Prof. F. Iavernaro \\ Dipartimento di Matematica, Universit\`a di Bari, via Orabona 4, 70125 Bari, Italy}

\abstract[Abstract]{We exploit the idea to use the maximal-entropy method, successfully tested in information theory and statistical thermodynamics, to determine approximating function's coefficients and squared errors' weights simultaneously as output of one single problem in least-square approximation. We provide evidence of the method's capabilities and performance through its application to representative test cases by working with polynomials as a first step. We conclude by formulating suggestions for future work to improve the version of the method we present in this paper.}


\jnlcitation{\cname{%
\author{Giordano D} and \author{Iavernaro F}} 
(\cyear{2020}), 
\ctitle{Maximal-entropy driven determination of weights in least-square approximation}, \cjournal{Mathematical Methods in Applied Sciences}, \cvol{2020;00:x--y}.}

\maketitle


\section{Introduction}\label{intro}

The method of least squares has a very interesting history that originated from the dispute for priority of discovery between the two giant mathematicians Legendre \cite{al1805} 
and Gauss \cite{cfg1809} in Europe and the (more or less) contemporary work of Adrain \cite{ra1808ta}, an Irish mathematician, in the United States. 
We enjoyed reading the neat historical accounts provided by Cleveland Abbe \cite{amca1871ajsa}, Merriman \cite{mm1877ta,mm1877}, Plackett \cite{rp1972b}, 
Stigler \cite{sms1977hm,sms1981as}, Dutka \cite{jd1990ahes}, and Aldrich \cite{ja1998isr}.
The basic idea of the method formulated by Legendre, Gauss and Adrain is conceptually rather simple and has been thoroughly elaborated upon and enriched during the course of the years \cite{mm1877,kp1901pmjs,mm1909,nc1920pmjs,rms1920pmjs,cgm2011,hsu1923josa,wed1930pps,wed1931pmjs,wed1934pmjs,wed1935pmjs,wed1948,gf1967jsiam,cwc1965jam,dy1966cjp,mp1967cj,mon1969cj,hm2006bmcb,ac2017jcm}.
One significant advance came with the introduction of the statistical weights.
Indeed, data may generally be affected by uncertainties arising from unpredictable, non-controllable, non-reproducible factors. 
Therefore, it seems not wise to treat all data within a given set as equivalent in the least-square minimization procedure; it makes more sense to accompany {\blue the square of each approximation error} 
with a multiplicative coefficient, the statistical weight, whose task is to gauge the 
relative importance {\blue of that particular square with respect to all others}.
But how are the statistical weights to be assigned?
This is a delicate matter whose importance was stressed long time ago by Uhler \cite[page 1065]{hsu1923josa}:
\begin{quotation} \noindent 
    If it be deemed worth while to use the method of least squares, one should first make a reasonable estimate of the relative weights of the observed quantities and then proceed to apply the method correctly.
\end{quotation}
and by Deming \cite[page 157]{wed1931pmjs}:
\begin{quotation} \noindent 
    The assignment of weights appears to be fundamental. Once the weights are assigned, the procedure is quite definitely marked out, and the results are unambiguous. 
    The question of just what values should be assigned to the weights in the first place is often not so clear cut. 
\end{quotation}
It is clear that the \textit{manual} assignment of the weights is nothing else than a subjective operation that introduces arbitrariness in the mathematical procedure and, therefore, should be avoided.
The literature offers several attempts \cite{mm1877,mm1909,wed1931pmjs,wed1934pmjs,hm2006bmcb,ac2017jcm} to remove that arbitrariness by exploiting reasoned manners to determine the weights within the algorithms designed to produce the approximating functions.
In the cited approaches, however, the weights are obtained from formulae based on approximated mathematical arguments.
In our approach, we take a different stand that draws inspiration from a pillar postulate of Shannon's mathematical theory of communication \cite{ces1948bstj3,ces1948bstj4} masterly phrased by Jaynes \cite[Sec. 3 on page 97]{etj1967} as follows:
\begin{quotation}  \noindent 
	If we accept \textsc{S\large{hannon}}'s interpretation (which can be justified by other mathematical arguments entirely independent of the ones given by \textsc{S\large{hannon}}) that the quantity
	\[  H = - \sum_{i} p_{i} \ln p_{i} \parbox{0em}{\hspace*{15em}  (8)}\]
	is an ``information measure'' for any probability distribution $p_{i}$, i.e. that it measures the ``amount of uncertainty'' as to the true value of $i$, then an ancient principle of wisdom --- that one ought to acknowledge frankly the full extent of his ignorance --- tells us that the distribution that maximizes $H$ subject to constraints which represents whatever information we have, provides the most honest description of what we know. 
	The probability is, by this process, ``spread out'' as widely as possible without contradicting the available information.
\end{quotation}
Thus, we reinterpret the weights as the probabilities in (8) of Jaynes' quote and proceed to determine the distribution of weights that maximizes the ``amount of uncertainty'' or entropy, in statistical-thermodynamics parlance, and subjected to two constraints: a) the weights are normalized and b) the minimal mean squared error is prescribed to a desired small value.
For the sake of completeness, we ought to mention that the concept of entropy maximization has already been used in connection with the least-square method but \textit{not} for the purpose of determining the statistical weights; in this regard, we refer the interested readers to the article of Preckell \cite{pp2001ajae} and the bibliography therein.
{\refereetwo Regarding a more global view about the applicability scope of the entropy-maximization method, besides the well known domains of information theory and statistical thermodynamics, we mention Wu's textbook \cite{MEP-mathpys-1995} for a theoretical analysis with applications in mathematics and physics and the interesting papers by Ravera \etal\ \cite{MEP-chem-2016} and by Mu\~{n}oz-Cobo \etal\ \cite{MEP-eng-2017} for recent applications in chemistry and engineering.}

\section{Theory}

Typically, we are given a set of data  $\{Y_{i},X_{i}\} \; (i=1,\ldots,n) \,$ and the problem we have to deal with is to find a function $\Fun{}$ such that the values $\Fun{i}$ approximate as better as possible the $Y_{i}$ data.
The approximating function can have polynomial, exponential, logarithmic, and so on, structure and, in general, includes an array $\bv{A}$ of $m+1$ unknown constant coefficients.
This formulation already implies either a simplification or a particular case: we either assume or are sure that there {\rosso are} no uncertainties affecting the $X_{i}$ data.
In this first introductory study, we proceed in this way for the sake of keeping contained the complicacies of the mathematical formalism and postpone the more general case to subsequent investigations. 
Also, rather than working directly with the given data, we prefer to preprocess them according to the transformation rules
\begin{subequations}\label{YX.transf}
  \begin{align}
     y_{i} & = \frac{Y_{i}-\Ymin}{\Ymax-\Ymin}  \label{YX.transf.y}   \\[.5\baselineskip]
     x_{i} & = \frac{X_{i}-\Xmin}{\Xmax-\Xmin}  \label{YX.transf.x} 
  \end{align}
\end{subequations}
to generate a set of adapted data $\{y_{i},x_{i}\} \; (i=1,\ldots,n) \,$; in \Req{YX.transf}, $\Ymin,\Xmin$ and $\Ymax,\Xmax$ are, respectively, the minima and the maxima in the given-data set.
We believe there are four good reasons that justify the preliminary adaptation of the given data. 
The first three are immediately evident from \Req{YX.transf}: 
a) the adapted data are non-dimensional numbers;
b) the adapted data are normalized because they are confined in the interval $[0,1]$; 
c) systematic errors, if any, affecting the given data are automatically swept out of the adapted data by the differences appearing on the right-hand sides.
The fourth reason is brought forward when we consider the appropriate definition of the {\rosso approximation error} involved in the replacement of $Y_{i}$ with $\Fun{i}$.
Intuitively, our attention is attracted by the difference $\Fun{i} - Y_{i}$ but that difference by itself does not say much whether or not the {\rosso error} brought about by replacing $Y_{i}$ with $\Fun{i}$ is small, and therefore good and acceptable, or great, and therefore bad and unacceptable; meaningfulness is acquired only when that difference is referred to a characteristic scale factor of the $Y_{i}$ data.
In this regard, we take the bandwidth $\Ymax-\Ymin$ as the most appropriate scale factor and adopt
\begin{equation}\label{err.def}
  e_{i} = \frac{\Fun{i}-Y_{i}}{\Ymax-\Ymin} 
\end{equation}
as operative definition for the {\rosso approximation error.} 
By taking into consideration \Req{YX.transf.y} and by introducing the scaled approximating function
\begin{equation}\label{sifun}
  \fun{} = \frac{\Fun{}-\Ymin}{\Ymax-\Ymin}
\end{equation}
the definition \Req{err.def} can be rephrased as 
\begin{equation}\label{err.def.1}
  e_{i} =  \fun{i} - y_{i}.
\end{equation}
On the contrary of the {difference} $\Fun{i} - Y_{i}$, the {\blue residual} \erase{difference} on the right-hand side of \Req{err.def.1} is extremely meaningful because, apart the absence of the irrelevant factor $100$, it gives the percent {\rosso error }implied by the {difference} $\Fun{i} - Y_{i}$ with respect to the bandwith $\Ymax-\Ymin$. 
That is the fourth reason supporting the convenience of working with the adapted-data set.
As a matter of fact, the target of our investigation shifts from $\Fun{}$ to $\fun{}$; the former, if needed, is provided by the inversion of \Req{sifun} after that the latter has been determined.

The core of the least-square method is the minimization of the weighted mean squared error
\begin{equation}\label{mse}
   \mse = \sum_{i=1}^{n} p_{i}\,\mse_{i} = \sum_{i=1}^{n} p_{i}\,[\fun{i} - y_{i}]^{2}
\end{equation}
that leads to the determination of the $m+1$ coefficients $\bv{a}$; the weights $p_{i}$ in \Req{mse} are subjected to the normalization condition
 \begin{equation}\label{wnc}
     \sum_{i=1}^{n} p_{i} = 1
 \end{equation}
The standard procedure {\blue consists in taking} the differential of \Req{mse} with respect to the coefficients $\bv{a}$, and {\blue setting} to zero the corresponding partial derivatives to generate an algebraic system
\begin{equation}\label{sys.a}
   \subeqn{\pd{}{\mse}{a_{k}} = 2 \sum_{i=1}^{n} p_{i}\,[\fun{i} - y_{i}] \pd{}{\fun{i}}{a_{k}} = 0}{k=0,\ldots,m}{}
\end{equation}
of $m+1$ equations for the unknown coefficients $\bv{a}$.
The weights behave as constant parameters during this series of operations.
The formal solution of \Req{sys.a} provides the coefficients
\begin{equation}\label{fs.a}
   \bv{a} = \bv{a}\,(p_{1},p_{2},\ldots,p_{n})
\end{equation}
and the minimal mean squared error
\begin{equation}\label{msem}
   \msem = \sum_{i=1}^{n} p_{i}\,[\funb{i} - y_{i}]^{2} {\blue = \msem\,(p_{1},p_{2},\ldots,p_{n})}
\end{equation}
 as functions of the weights.
 This would be the finishing line if the weights were assigned.
 But we assume they are not and, in order to find them, we extend the least-square standard procedure in accordance with the principle stated in Jaynes' quote at the end of \Rse{intro}: we look for the weight distribution {\blue compatible with the constraints \Req{wnc} and \Req{msem}} that maximizes the ``amount of uncertainty'' or entropy
 \begin{equation}\label{aou}
    H = - \sum_{i=1}^{n} p_{i} \ln p_{i}
 \end{equation} 
Within this perspective, the minimal mean squared error $\msem$ {\blue in \Req{msem}} must be supposed as prescribed.
Thus, we build the Lagrangian
\begin{equation}\label{lag}
   L = \sum_{i=1}^{n} p_{i} \ln p_{i} + \alpha \cdot \left( \sum_{i=1}^{n} p_{i} - 1 \right) + \beta \cdot \left\{ \sum_{i=1}^{n} p_{i}\,[\funb{i} - y_{i}]^{2} - \msem    \right\}
\end{equation}
and we take its differential with respect to the weights by paying utmost attention to the functional dependence \Req{fs.a}
\begin{equation}\label{dlag}
   dL = \sum_{i=1}^{n} (1 + \ln p_{i})dp_{i} 
      + \alpha \sum_{i=1}^{n} dp_{i} 
      + \beta  \sum_{i=1}^{n} [\funb{i} - y_{i}]^{2} dp_{i}
      +2\beta  \sum_{i=1}^{n} \left\{ p_{i}\,[\funb{i} - y_{i}] \sum_{j=1}^{n} dp_{j} \sum_{k=0}^{m} \pd{}{a_{k}}{p_{j}}\,\pd{}{\funb{i}}{a_{k}} \right\}      
\end{equation}
Fortunately, the last term on the right-hand side of \Req{dlag} vanishes identically
\begin{equation}\label{dlag.last}
   \sum_{i=1}^{n} \left\{ p_{i}\,[\funb{i} - y_{i}] \sum_{j=1}^{n} dp_{j} \sum_{k=0}^{m} \pd{}{a_{k}}{p_{j}}\,\pd{}{\funb{i}}{a_{k}} \right\} =
   \sum_{j=1}^{n} dp_{j} \left\{ \sum_{k=0}^{m} \pd{}{a_{k}}{p_{j}} \sum_{i=1}^{n} p_{i}\,[\funb{i} - y_{i}]  \pd{}{\funb{i}}{a_{k}} \right\} = 0
\end{equation}
because of \Req{sys.a}.
After this simplification, the {\blue differential \Req{dlag} of the Lagrangian} reduces to
\begin{equation}\label{dlag.s}
   dL = \sum_{i=1}^{n} \left\{ 1 + \ln p_{i} + \alpha + \beta \, [\funb{i} - y_{i}]^{2} \right\} dp_{i} 
\end{equation}
from which we can extract and set to zero the partial derivatives
\begin{equation}\label{sys.p}
   \subeqn{\pd{}{L}{p_{i}} = 1 + \ln p_{i} + \alpha + \beta \, [\funb{i} - y_{i}]^{2} = 0}{i=1,\ldots,n}{}
\end{equation}
to generate an algebraic system of $n$ equations for the unknown weights.
We have now available a coupled set of $m+1$ equations \Req{sys.a}, $n$ equations \Req{sys.p} and two constraints \Req{wnc}, \Req{msem} for the determination of $m+1$ coefficients $\bv{a}$, $n$ weights $p_{i}$ and two Lagrangian multipliers $\alpha, \beta$.
In this way, the seemingly separate necessities of finding the coefficients of the approximating function and of assigning the weights merge into one single mathematical problem.
{\refereeone
Before continuing, we believe that an important consideration is in order here.
The formula (\ref{aou}) suggests that the maximal-entropy procedure we have developed is nonlinear in nature even in the event that the approximating function $\fun{}$ depends linearly on the coefficients $\bv{a}$. 
This aspect has a relevant mathematical implication. 
Assuming prescribed positive weights and excluding some special degenerate cases, it is well known that, in a standard linear least-square problem, the mean squared error $\msem$ in (\ref{msem}) is a strictly convex function of the coefficients $\bv{a}$ and its minimum is a unique solution of the normal linear system \Req{sys.a}. 
If the maximal-entropy idea is exploited to determine the weights, the corresponding function to optimize becomes, in general, non-convex; this occurrence leads to the existence of multiple optimal solutions. 
The necessity of a deeper theoretical study of this matter being understood, the algorithmic procedure described at the end of the present section [paragraph beginning after \Req{msem.uw}] clarifies the extent to which the optimal solution we are seeking improves on the one associated with the classical least-square procedure with constant weights.
}

The constant $1+\alpha$ involving the first Lagrangian multiplier is eliminable by first resolving \Req{sys.p} for $p_{i}$
and then imposing the normalization condition \Req{wnc}; this two-step operation leads to the definition of the partition function\footnote{We are deliberately borrowing the term from the statistical-thermodynamics parlance.}
\begin{equation}\label{pf}
   Q = \sum_{i=1}^{n} \exp\left\{ - \beta \, [\funb{i} - y_{i}]^{2} \right\}
\end{equation}
to the explicit expression of the constant $1+\alpha$
\begin{equation}\label{alpha}
   \exp(1 + \alpha) = Q
\end{equation}
and to the final form of the algebraic system \Req{sys.p}
\begin{equation}\label{sys.p.f}
   \subeqn{p_{i} = \frac{1}{Q} \exp\left\{ - \beta \, [\funb{i} - y_{i}]^{2} \right\}}{i=1,\ldots,n}{}
\end{equation}
We believe it is important here to stress that the functional dependence \Req{fs.a} is implied on the right-hand sides of both \Req{pf} and \Req{sys.p.f}. 

The second Lagrangian multiplier $\beta$ is not eliminable, unfortunately.
The equation needed for its determination is obtained by substituting {\blue \Req{sys.p.f}} into \Req{msem}, with due account of \Req{pf}, and {\blue by} rearranging to the final form
\begin{equation}\label{beta}
   \msem \cdot \sum_{i=1}^{n} \exp\left\{ - \beta \, [\funb{i} - y_{i}]^{2} \right\} = 
   \sum_{i=1}^{n} [\funb{i} - y_{i}]^{2} \cdot \exp\left\{ - \beta \, [\funb{i} - y_{i}]^{2} \right\}
\end{equation}
From a numerical-calculation point of view, \Req{beta} is handled extremely well by the Newton-Raphson iterative method which converges very quickly. 

There is one more task to take care of: to express the entropy \Req{aou} as a function of the mean squared error.
Substitution of \Req{sys.p.f} into \Req{aou} and further manipulation with due account of \Req{wnc} and \Req{beta} leads to the relatively simple expression
\begin{equation}\label{Hexpl}
   H = H(\msem) = \beta \cdot \msem + \ln Q
\end{equation}
from which, in turn, we obtain by derivation the rather important result
\begin{equation}\label{dHexpl}
   \tds{}{H}{\msem} = \beta 
\end{equation}
We skip the derivation details to save space and limit ourselves to warn the reader who wishes to try out the enjoyment of the exercise to pay extreme attention when carrying out the subtle derivation of $\ln Q$; \Req{sys.a} will again come to help as it did in \Req{dlag.last}. 

For the readers to whom the maximal-entropy nuance we have introduced in the least-square method may appear somewhat unconvincing, we wish to show that our path does not diverge from the one traced by the standard paradigm but that we simply go one step further along the same direction.
If we set $\beta=0$ in \Req{lag} then we basically relinquish, or better ignore, the constraint \Req{msem} and content ourselves only with the normalization condition \Req{wnc}; in such a case, the partition function \Req{pf} reduces to
\begin{equation}\label{pf.uw}
   Q = n
\end{equation}
and the algebraic system \Req{sys.p.f} yields the simple uniform-weight solution 
\begin{equation}\label{sys.p.uw}
   \subeqn{p_{i} = \frac{1}{n}}{i=1,\ldots,n}{0}
\end{equation}
The {\blue uniform-weight} coefficients $\bv{a}_{uw}$ follow from \Req{sys.a} after substitution of \Req{sys.p.uw}; the minimal mean squared error from \Req{msem} reduces to the standard average
\begin{equation}\label{msem.uw}
   \msem_{uw} = \frac{1}{n}\sum_{i=1}^{n} [\funuw{i} - y_{i}]^{2}
\end{equation}
Thus, we have retrieved the uniform-weight solution of the least-square method and, at the same time, have proven that also that solution is a consequence of entropy maximization.
Our inclusion of the constraint \Req{msem} represents, therefore, just an additional enhancement, still compliant with the same guideline leading to the uniform-weight solution.
After having pointed out this clarification, from hereinafter we will use the term maximal-entropy method, and acronym MEM, to refer exclusively to our approach based on the inclusion of \Req{msem} as constraint of entropy maximization to obtain simultaneously coefficients $\bv{a}$ and weights $p_{i}$.

To summarize, the $m+1$ equations \Req{sys.a}, the $n$ equations \Req{sys.p.f} and their auxiliary equation \Req{pf}, and the constraint-derived equation \Req{beta} 
provide, respectively, the coefficients $\bv{a}$, the weights $p_{i}$, and the Lagrange multiplier $\beta$; this is the core of MEM.
The number of equations equals the number of unknowns and we are set to switch on the numerical-calculation machinery.
Unfortunately, the set of equations is highly nonlinear.
The solution strategy we have devised consists of the following steps:
\begin{enumerate}
   \item   Assign 
           \begin{itemize}
              \item[a)]   an initial weight distribution, likely the uniform-weight distribution
              \item[b)]   the desired minimal mean squared error $\msem$
           \end{itemize}   
   \item   \label{two} Evaluate the coefficients $\bv{a}$ from \Req{sys.a}
   \item   Update the Lagrangian multiplier $\beta$ from  \Req{beta}
   \item   Update the weights $p_{i}$ from \Req{pf} and \Req{sys.p.f}
   \item   Repeat from step \ref{two} until convergence is achieved.
\end{enumerate}
It is interesting to notice that, when convergence is reached, \Req{sys.p.f} assures the Gaussian nature of the weight distribution with respect to the {\rosso errors} \Req{err.def.1}.
Moreover, a clarification is in order regarding step {\rosso 1.b)}.
We have seen already [paragraph beginning after \Req{dHexpl}] that the uniform-weight distribution is that particular solution which implies $\beta=0$.
This is an important distribution and it should be considered as a reference baseline.
The value $\msem_{uw}$ returned by \Req{msem.uw} constitutes the reference level with respect to which the desired minimal mean squared error, mentioned in step {\rosso 1.b)} and playing a role in \Req{beta}, should be assigned, that is, $\msem < \msem_{uw}$.
As a matter of fact, in order to keep this inequality well under control in the implementation of our solution strategy, we prefer to assign a reduction factor $r>1$ and to set $\msem = \msem_{uw} / r$, rather than assigning directly $\msem$.

We have applied systematically our solution strategy to several test cases by taking a polynomial approximating function 
\begin{equation}\label{ifpol}
   \fun{} = \sum_{k=0}^{m} a_{k} x^{k}
\end{equation}
to begin with, because the derivatives
\begin{equation}\label{der.a}
   \pd{}{\fun{}}{a_{k}} = \sum_{s=0}^{m} \pd{}{a_{s}}{a_{k}} x^{s} = x^{k}
\end{equation}
do not depend on the coefficients and the algebraic system \Req{sys.a} required for their calculation becomes linear 
\begin{equation}\label{sys.a.pol}
   \subeqn{ \sum_{s=0}^{m}\left(\sum_{i=1}^{n} p_{i}\,x_{i}^{k} x_{i}^{s}\right) a_{s} = \sum_{i=1}^{n} p_{i}\,y_{i} x_{i}^{k}}{k=0,\ldots,m}{}
\end{equation}
We have obtained satisfactory results that helped us to acquire good confidence about the solution strategy's performance and, more in general, in the maximal-entropy prescript. 
In \Rse{ntc}, we will describe four test cases.
The first one is rather simple and is meant to put in evidence differences and advantages of the maximal-entropy enhanced procedure with respect to the uniform-weight standard procedure.
The other three have been selected to show the usefulness of the maximal-entropy procedure from the perspective of potential application in different sectors.
All our results have been obtained by a numerical algorithm  ({\rosso labeled} FM in the forthcoming figures, where needed) based on the five-step strategy described above  and programmed in \texttt{Octave}.

\section{Numerical test cases}\label{ntc} 

\subsection{Pearson's linear-approximation problem}\label{kptc}

Pearson's linear-approximation problem \cite[bottom of page 569]{kp1901pmjs} is rather simple and, because of its simplicity, it has been repeatedly considered in the literature as benchmark problem. 
It consists in finding ``the best fitting straight line'' to the data shown in the two left-most columns of \Rta{kpd}. 
Pearson considered three cases categorized according to the existence of uncertainties, respectively, in both $Y_{i}$ and $X_{i}$ data, only in the $Y_{i}$ data, and only in the $X_{i}$ data. 
For compliance with our initial assumption of no uncertainties affecting the $X_{i}$ data, therefore, we concentrated our analysis only on the second case; the numerical details can be found in \cite[page 570]{kp1901pmjs}.
\begin{table}[h]
   \caption{Pearson's data \cite[bottom of page 569]{kp1901pmjs}.\label{kpd}}  \centering
   \resizebox{.45\textwidth}{!}{\includegraphics*{\figdir/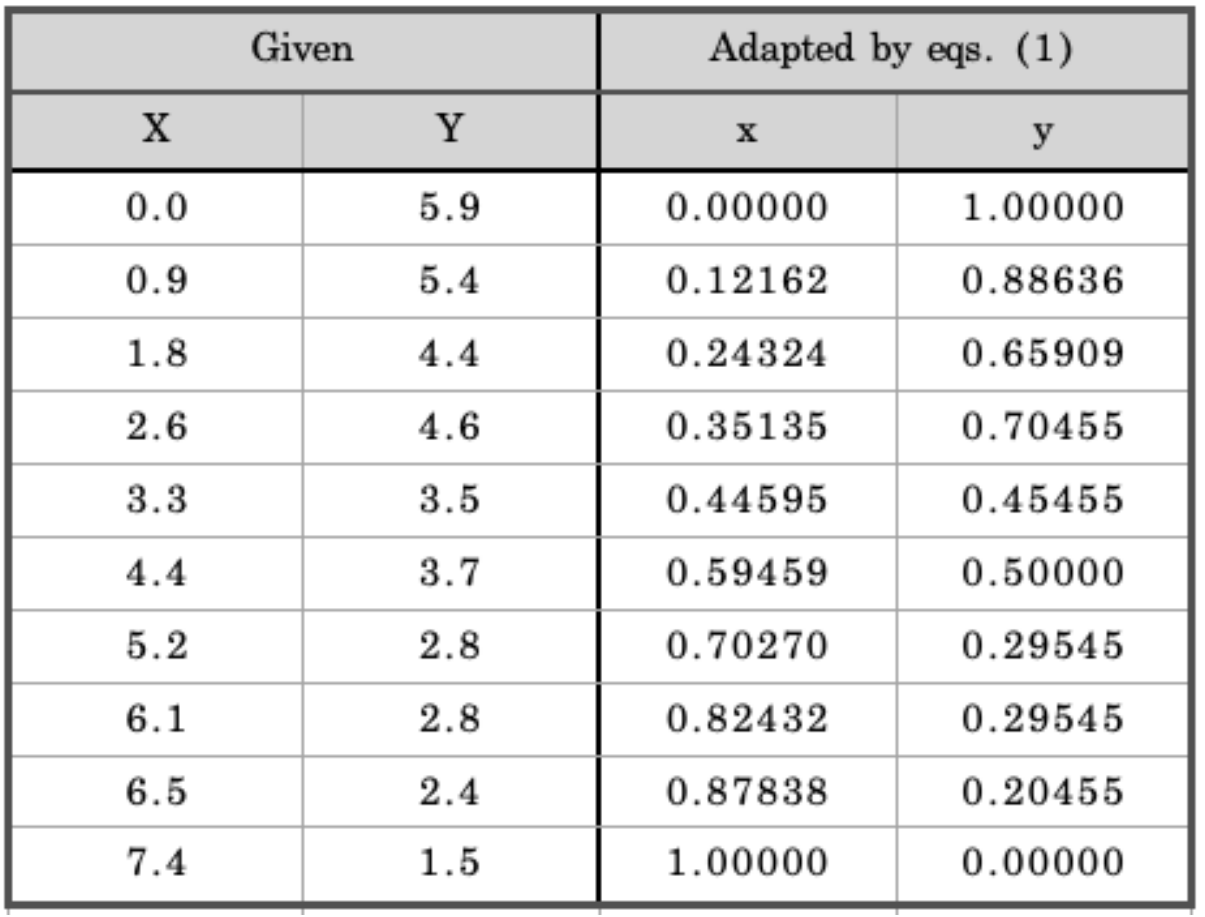}}
\end{table}
\begin{figure}[h]
   \centering
   \resizebox{.500\textwidth}{!}{\includegraphics*[trim=17 25 50 50]{\figdir/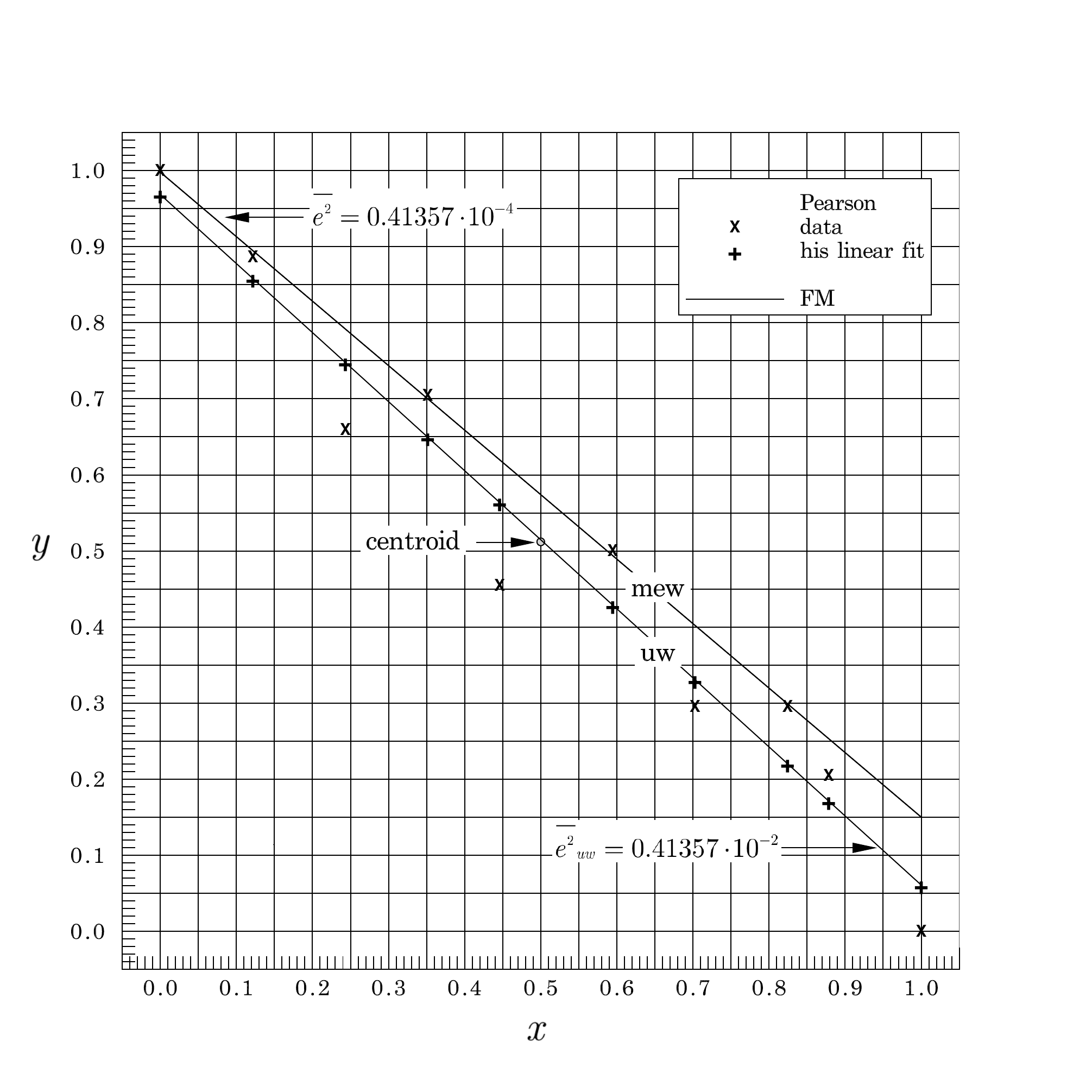}}  \\
   \resizebox{.495\textwidth}{!}{\includegraphics*[trim=17 15 55 50]{\figdir/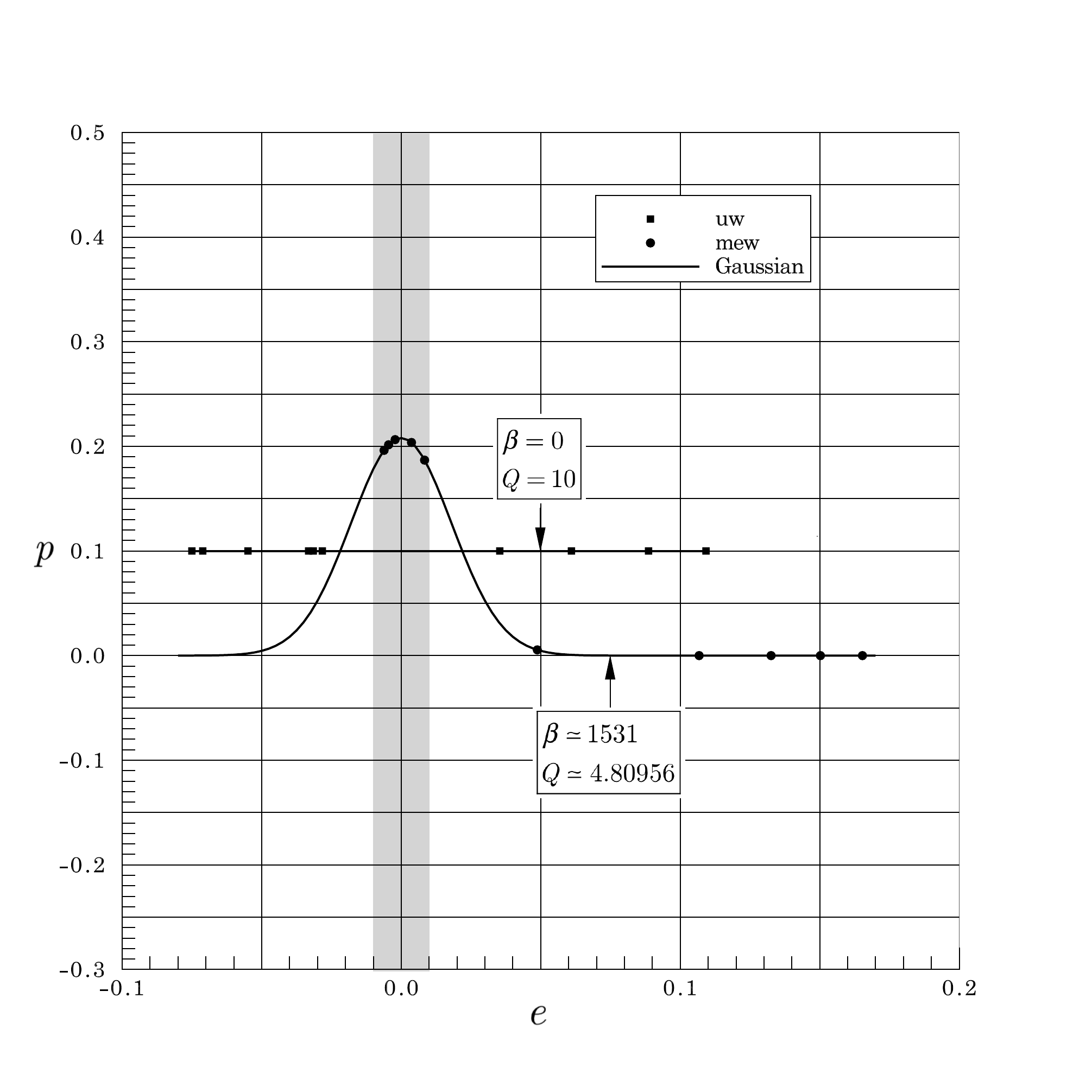}} \hfill
   \resizebox{.495\textwidth}{!}{\includegraphics*[trim=17 15 55 50]{\figdir/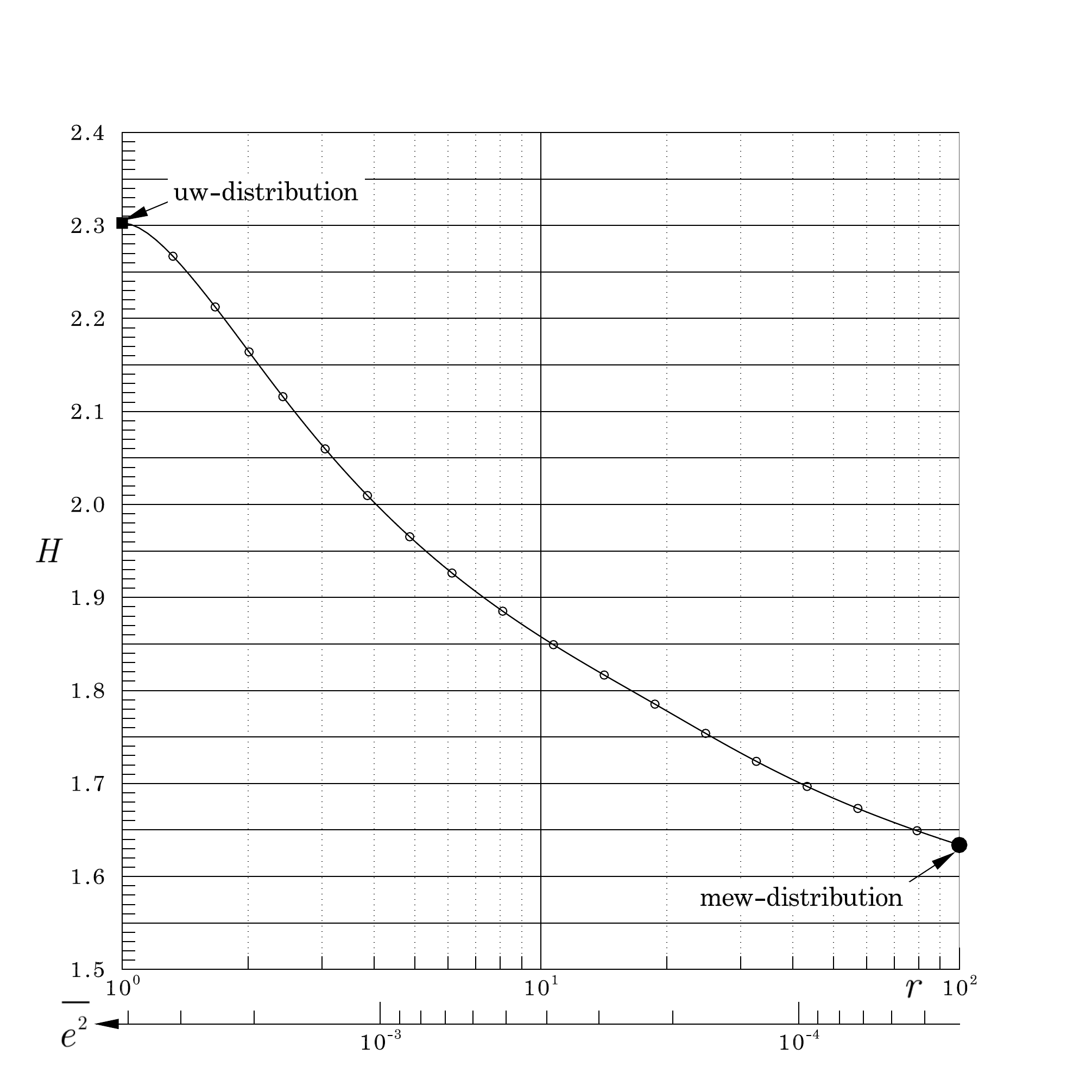}}
   \caption{Top:   straight-line approximation of Pearson's data \cite[figure at page 570]{kp1901pmjs} according to uniform and maximal-entropy weights.
            Left:  weight distributions \Req{sys.p.f} versus {\rosso residuals} \Req{err.def.1} corresponding to the uw and mew lines of top diagram.
            Right: entropy history as a function of either reduction factor or mean squared error during transition from uw- to mew-distribution.\label{kpwd}}
\end{figure}
The comparison between Pearson's and our results are illustrated in \Rfi{kpwd} (top), styled to some extent similarly to Pearson's original figure (see \cite[page 570]{kp1901pmjs}). 
Of course, there is full coincidence between his linear fit (``\scalebox{0.8}{\textbf{+}}'' symbols) and our uniform-weight (uw) line
\begin{equation}\label{p-uw}
   f(x)_{uw} = -0.90747 x + 0.96845
\end{equation}
to which there corresponds a minimal mean squared error $\msem_{uw}= 0.41357\cdot 10^{-2}$.
Yet, the uw-line does not go over any of the \erase{given} {\blue adapted} data (``\scalebox{0.8}{\textsf{X}}'' symbols) and its acceptance as \erase{interpolating} {\blue approximating} function implies substantial {\rosso  errors}; indeed, the worst of them occurs in correspondence to the fifth (from left) data point and turns out to be almost 11\%, perhaps not so convincingly negligible as wished.
This is the toll to pay for giving same importance, i.e. same weight, to all data.
MEM tells a completely different story.
By decreasing the desired minimal mean squared error from the initial $\msem_{uw}= 0.41357\cdot 10^{-2}$ to the final $\msem = 0.41357\cdot 10^{-4}$ ($r=100$), the approximating line will transition from the uw-line to the maximal-entropy-weight (mew) line 
\begin{equation}\label{p-mew}
   f(x)_{mew} = -0.84758 x + 0.99781
\end{equation}
The striking characteristic of the mew-line is the ability to go over five data points with {\rosso errors} less than 1\%. 
The presence of the other five data points would disturb the alignment and, therefore, the entropy maximization constrained by \Req{msem} filters them out by making their weights vanish.
This conjecture finds clear confirmation in the diagram of the weight distributions \Req{sys.p.f} versus {\rosso residuals} \Req{err.def.1} shown in \Rfi{kpwd} (left).
The mew-distribution (solid circles) presents all the five data points aligned on the mew-line in \Rfi{kpwd} (top) positioned within the narrow gray band [-1\%,\raisebox{0.25ex}{\scalebox{0.6}{+}}1\%], each one with its own weight; the bad approximations of the other five data points do not matter because they are shut off by vanishing weights.
The uw-distribution (solid squares) presents a larger approximation spread \mbox{[-7.5\%,11\%]}; of course the weights are all positioned on the horizontal line $p=1/n=1/10=0.1$ which, by the way, can be interpreted as a (degenerate) Gaussian distribution.
The diagram in \Rfi{kpwd} (right) illustrates the entropy history versus either reduction factor or mean squared error during the transition from the uw-distribution to the mew-distribution.

The ability of MEM to filter out undesired data raised our curiosity and inspired a numerical experiment meant to check whether or not such an ability was just casual for Pearson's test case or, as a matter of fact, something more intrinsic to the method.
We considered the line 
\begin{equation}\label{y=x}
   y = x
\end{equation}
with $x$ in the interval $[0,1]$ and discretized the latter in $n = 21$ equidistant points.
We scrambled randomly the $y_{i}$ values of the even-positioned ($i=2,4,6,\ldots$) points and left untouched the odd-positioned points.
\begin{figure}[h]
   \centering
   \resizebox{.500\textwidth}{!}{\includegraphics*[trim=17 25 50 50]{\figdir/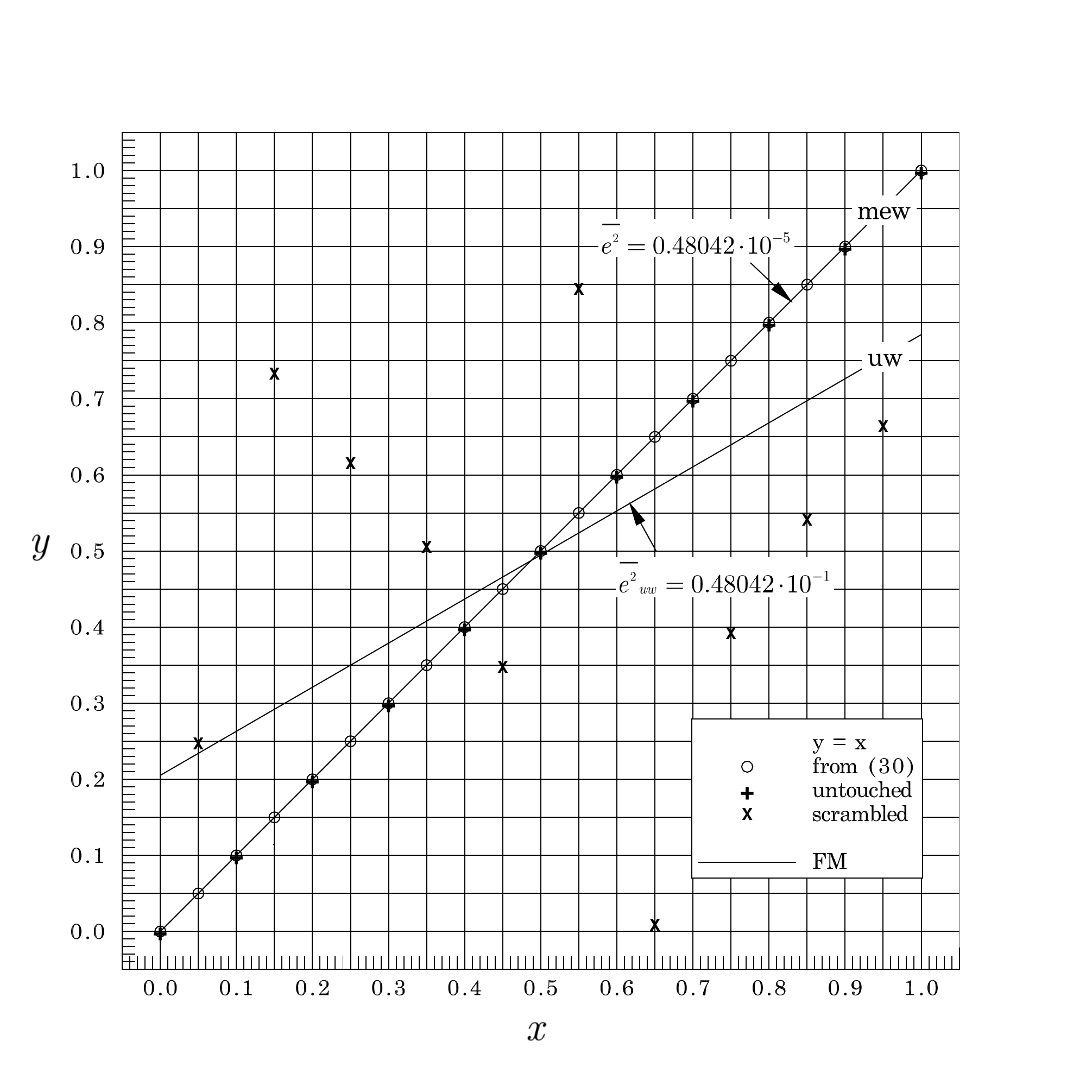}}  \\
   \resizebox{.495\textwidth}{!}{\includegraphics*[trim=17 15 55 50]{\figdir/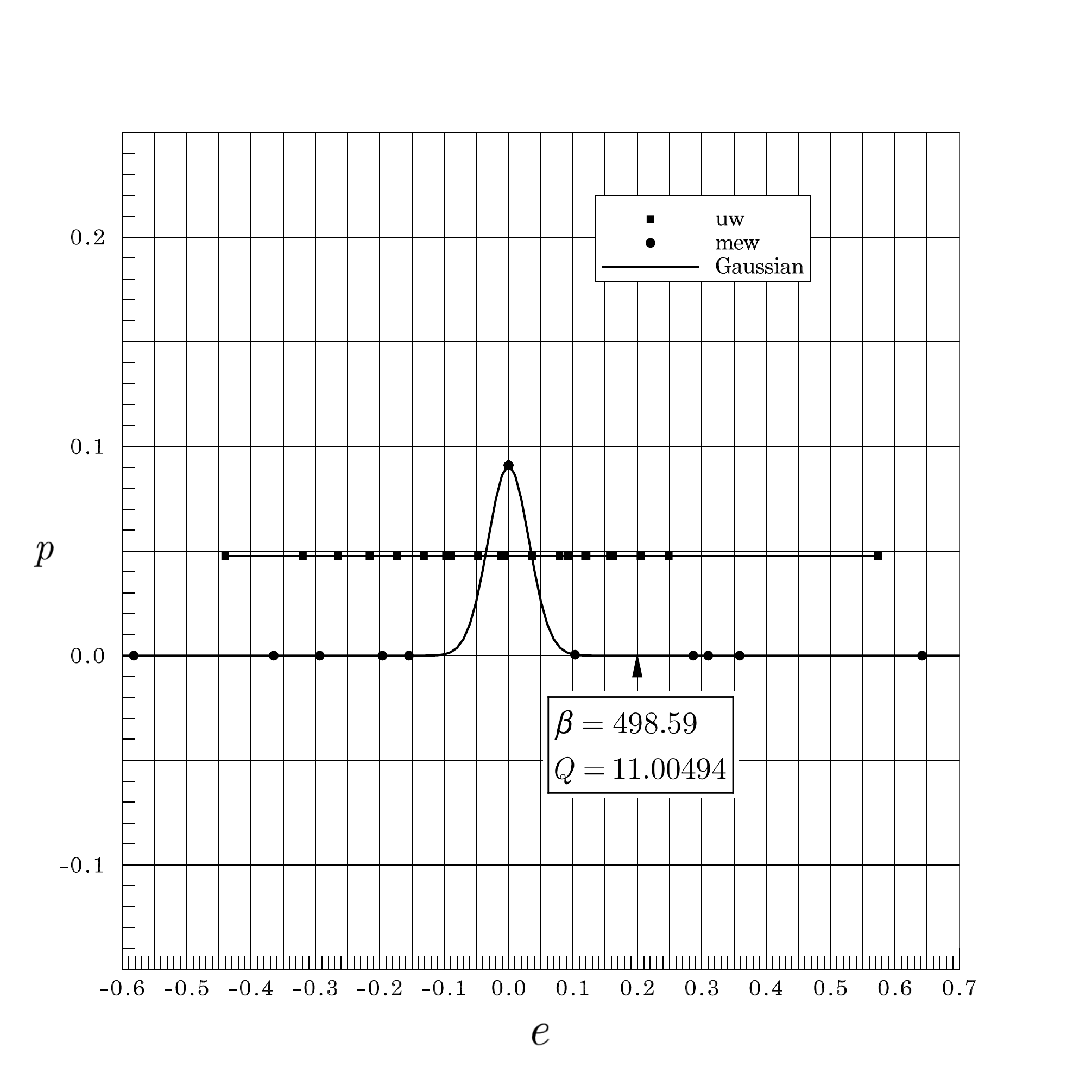}} \hfill
   \resizebox{.495\textwidth}{!}{\includegraphics*[trim=17 15 55 50]{\figdir/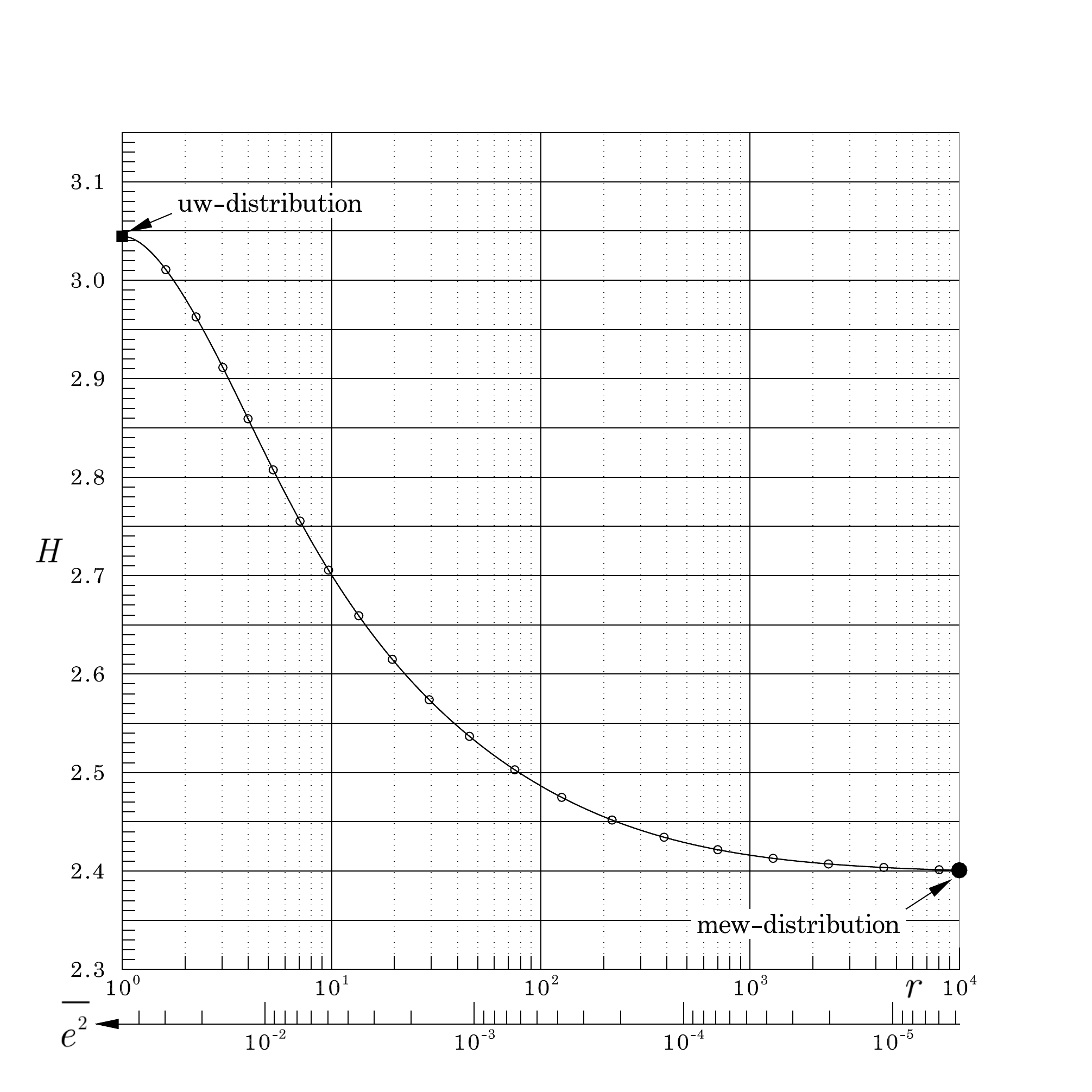}}
   \caption{Top:   least-square approximation of the hidden-line numerical experiment according to uniform and maximal-entropy weights.
            Left:  weight distributions \Req{sys.p.f} versus {\rosso residuals} \Req{err.def.1} corresponding to the uw and mew lines of top diagram. 
            Right: entropy history as a function of either reduction factor or mean squared error during transition from uw- to mew-distribution.\label{f1wd}}
\end{figure}
The results from our \texttt{Octave} program are shown in \Rfi{f1wd} (top).
The uw-line
\begin{equation}\label{f1-uw}
   f(x)_{uw} = 0.57920 x + 0.20512
\end{equation}
does, of course, its best to embrace both scrambled and untouched data points but its score $\msem_{uw}= 0.48042\cdot 10^{-1}$ is maybe not so good.
On the other hand, the mew-line 
\begin{equation}\label{f1-mew}
   f(x)_{mew} = 1.0000 x - 0.57865\cdot10^{-4}
\end{equation}
ignores the scrambled data and aligns diligently with the untouched ones by scoring, in so doing, an impressive $\msem = 0.48042\cdot 10^{-5}$.
The weight distributions versus {\rosso residuals \erase{(was approximation errors)}} are shown in \Rfi{f1wd} (left).
We see clearly that the scrambled data points belonging to the mew-distribution (solid circles) are all sitting on the horizontal line $p\simeq0$. 
The untouched data points are all coinciding with the solid circle sitting on the peak of the Gaussian curve; there are very slight differences in weight and {\rosso  residual \erase{(was error)}} among them but they are so minute that do not show in the scale of the figure.
The entropy history during transition, shown in \Rfi{f1wd} (right), displays an interesting asymptotic behavior for $r\rightarrow\infty$ or $\msem\rightarrow 0$: it levels at $H = \ln 11 \simeq 2.3979$ which is the value the entropy would attain if only the 11 untouched data points were available for interpolation and, logically, considered with equal weight ($p=1/11$) because they are perfectly aligned.
In conclusion, the outcome of our numerical experiment reinforces the conviction that the ability to filter out undesired data is a feature that belongs to MEM.

\subsection{Signal rebuilding}
Encouraged by the positive outcome of the hidden-line experiment described in the last part of \Rse{kptc}, we decided to add more complexity to the interpolation process:
we aimed at checking whether or not MEM could rebuild a signal partially contaminated by white noise. 
In honor of one of the discoverers of the least-square method, we assumed as signal the Legendre polynomial of degree seven
\begin{equation}\label{Lp7} \rosso
   y = \Lpsev = 1716 x^{7} - 6006 x^{6} + 8316 x^{5} - 5775 x^{4} + 2100 x^{3} - 378 x^{2} + 28 x
\end{equation}
purposely shifted to the interval $[0,1]$. 
Then, we discretized the latter with a uniform mesh
\begin{equation}\label{lp7m}
   \subeqn{x_{i}=\frac{i - 1}{99}}{i=1,\dots,100}{}
\end{equation}
composed by 100 equidistant points and randomly extracted a subset $\Omega$, with cardinality $\mbox{card}(\Omega)=75$, from the set of indices $(1,2,\dots 100)$ to produce the noise-contaminated $y_{i}$ values according to the selection rule
\begin{figure}[h]
   \centering
   \resizebox{.500\textwidth}{!}{\includegraphics*[trim=17 25 50 50]{\figdir/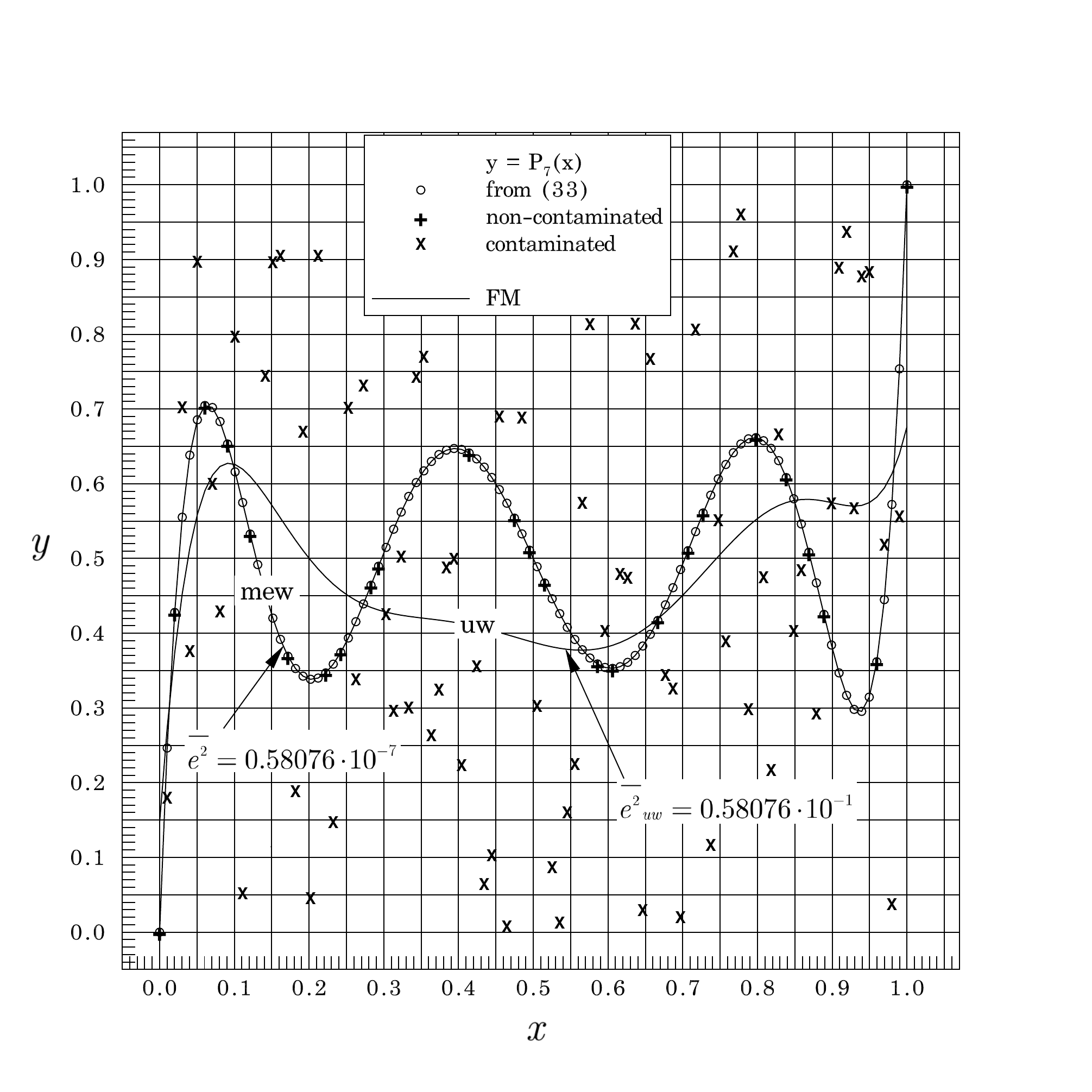}}
   \caption{Least-square approximation of the Legendre-polynomial numerical experiment according to uniform and maximal-entropy weights.\label{f3df}}
\end{figure}
\begin{equation}
   y_{i} = 
	   \begin{cases} \; P_7(x_i)      & \quad\text{if $i\not\in\Omega$}  \\[1.5ex] 
	                 \; \delta_{i}    & \quad\text{if $i    \in\Omega$}  \\[0ex] 
	   \end{cases} \label{perturbed-signal}
\end{equation}
In \Req{perturbed-signal}, $\delta_i$ represent random numbers uniformly distributed on the interval $[0, 1]$. 
\Rfi{f3df} illustrates the {\blue results: the hollow circles indicate the values calculated directly from \Req{Lp7}, the ``\scalebox{0.8}{\textbf{+}}'' symbols single out the 25 values not contaminated by noise (they coincide with the hollow circles, obviously), the ``\scalebox{0.8}{\textsf{X}}'' symbols refer to the 75 noise-contaminated values, the solid lines correspond to the output from our FM algorithm for both ``uw'' and ``mew'' cases.}
Once again, the uw-curve 
\begin{equation}\label{f3-uw}
   f(x)_{uw} = 375.28004 x^{7} - 1379.7699 x^{6} + 2036.3681 x^{5} - 1545.1538 x^{4} + 640.82981 x^{3} - 140.84931 x^{2} + 13.824592x +  0.14675605
\end{equation}
struggles hard amidst the jungle of the noise-contaminated data points and scores poorly with $\msem_{uw} = 5.8076\cdot 10^{-2}$.
In turn, the mew-curve 
\begin{equation}\label{f3-mew}
f(x)_{mew} = 1716.0167 x^{7} - 6006.0547 x^{6} + 8316.0673 x^{5} - 5775.0373 x^{4} + 2100.0078 x^{3} - 377.99964 x^{2} + 27.999760 x - 0.0000037
\end{equation}
corresponding to a reduction factor $r=10^{6}$, which implies $\msem = 5.8076\cdot 10^{-8}$, settles right on target. 
Visual {\blue inspection of both} \Req{Lp7} and \Req{f3-mew} 
reveals the extreme closeness of the corresponding polynomial coefficients and corroborates the satisfactory performance of the mew-curve {\blue displayed in \Rfi{f3df}}.

\subsection{Outlier removal}\label{subsec:out-rem}

A general working assumption in least-square approximation is that data are affected by errors which belong to a normal distribution. 
\begin{figure}[h]
   \centering
   \resizebox{.495\textwidth}{!}{\includegraphics*[trim=17 15 55 50]{\figdir/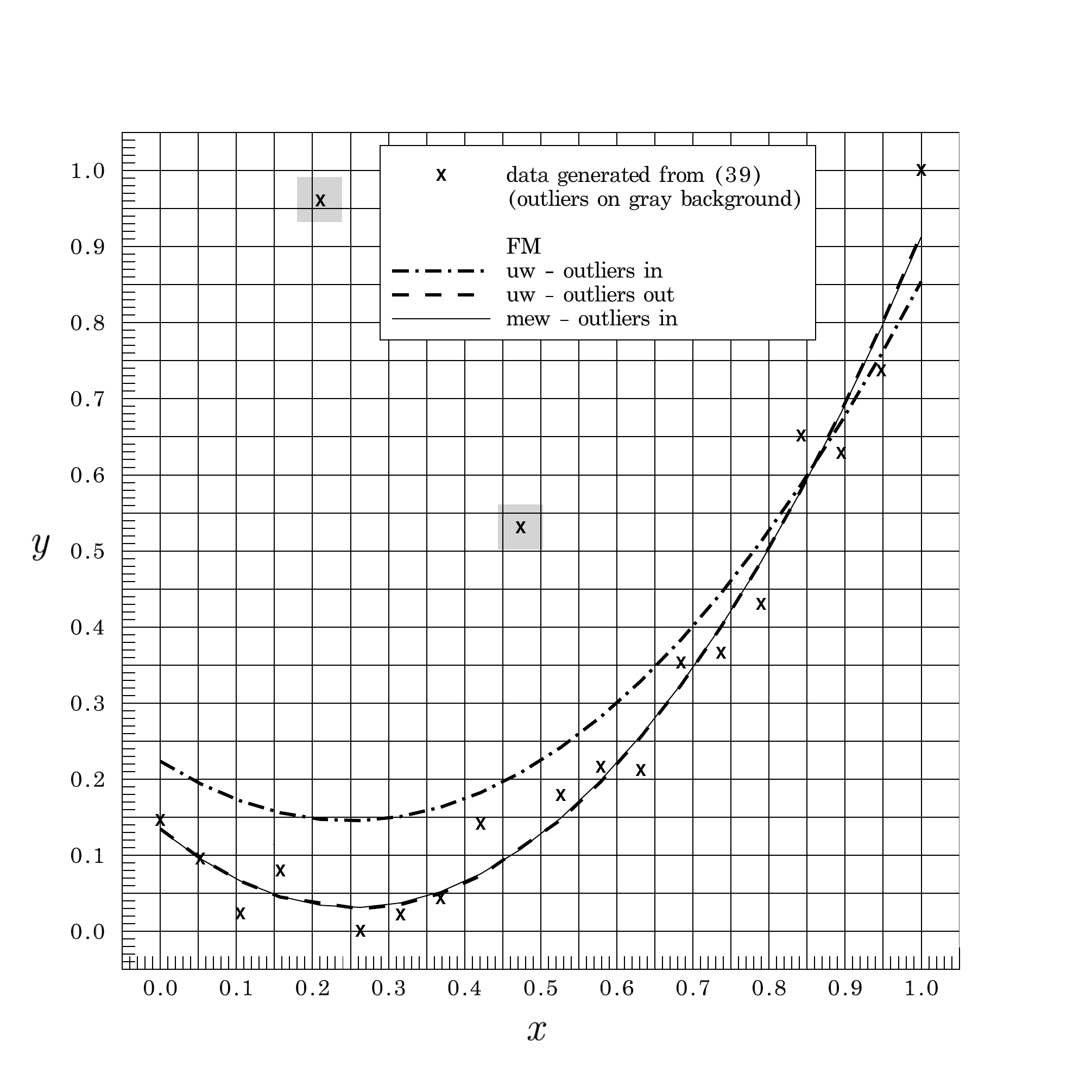}} \hfill
   \resizebox{.495\textwidth}{!}{\includegraphics*[trim=17 15 55 50]{\figdir/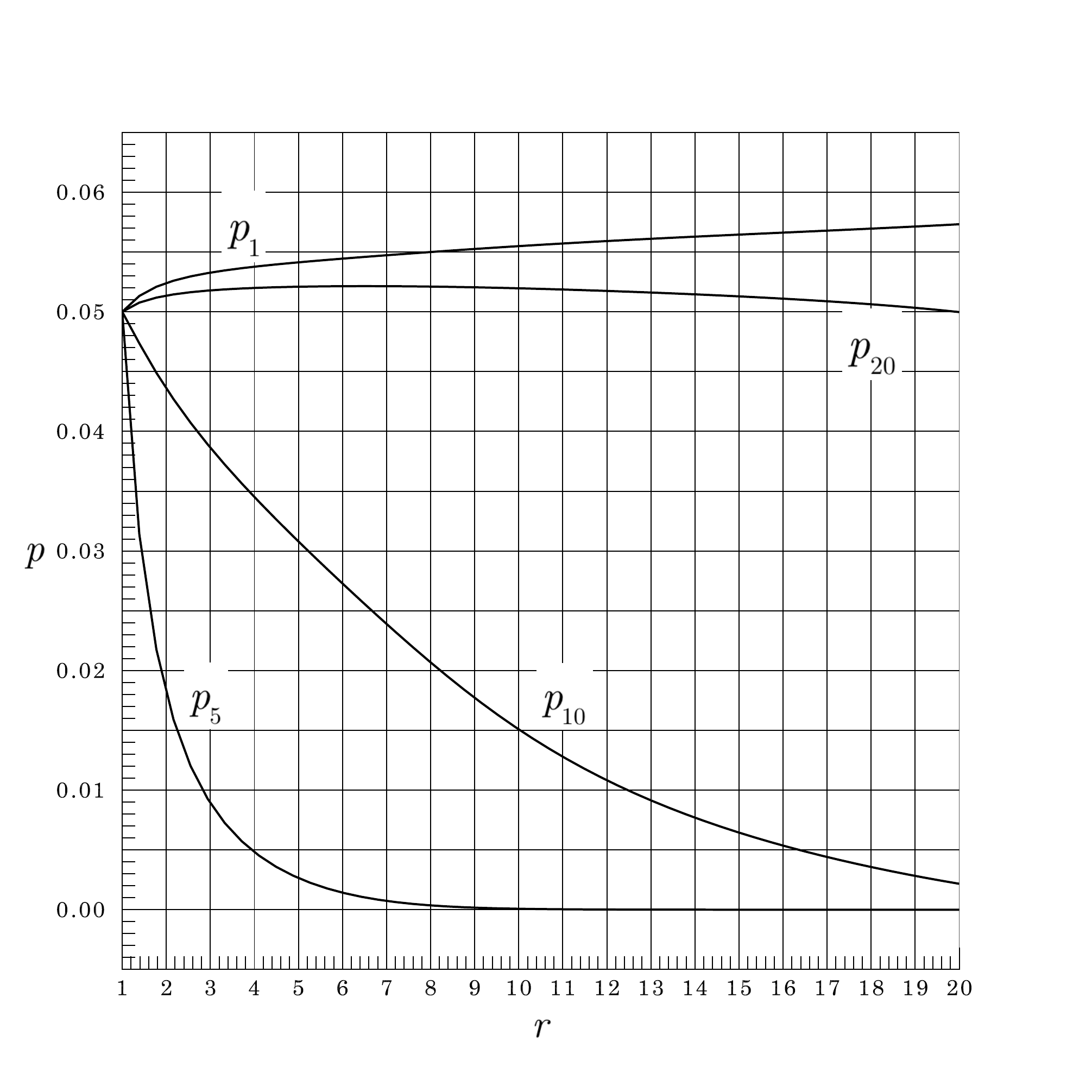}}
   \caption{Left:   least-square approximation of the perturbed parabola according to uniform (outliers included and manually removed) and maximal-entropy weights.
            Right:  maximal-entropy weights versus reduction factor.}
   \label{pic-outliers}
\end{figure}
However, the presence of outliers cannot always be ruled out; even a few of them may have a relevant impact on the shape of the approximating curve and may lead to misleading fitting \cite{hm2006bmcb,HMT83,HRRS06}. 
Outlier identification and removal are, therefore, necessary and critical operations.
Inspired by the interesting article by Motulsky and Brown \cite{hm2006bmcb}, we conceived the test case described in this section to show the ability of MEM to execute those tasks. 
We consider the parabola 
\begin{equation}\label{para}
   Y = 1 - X + 2 X^2
\end{equation}
with $X$ in the interval $[0,1]$. 
Then we discretize the latter in $n = 20$ equidistant points, introduce the perturbation
\begin{equation}\label{perturbed-data}
   \subeqn{Y_{i} = 
	   \begin{cases} \; Y(X_i) + \dfrac{\delta_{i}}{20}  & \quad \text{if $i\not=5,10$}  \\[1.5ex] 
	                 \; 2                                & \quad \text{if $i=5$}         \\[1.5ex] 
	                 \; 1.5                              & \quad \text{if $i=10$}        \\[0ex] 
	   \end{cases}}{i=1,\ldots,20}{} 
\end{equation}
with $\delta_i \sim {\cal N}(0,1)$ being a random variable distributed normally with mean $0$ and variance $1$, and finally feed the $\{Y_{i},X_{i}\} \; (i=1,\ldots,n) \,$ data into the transformation \Req{YX.transf} to produce the data in our working format.
The results are shown in \Rfi{pic-outliers}. 
On the left diagram, the dash-dot uw-parabola expectedly misknows the outliers $y_5$ and $y_{10}$ (on gray background) and is clearly affected by their presence; its mean squared error is $\msem_{uw}=0.46274\cdot 10^{-1}$. 
The mew-parabola (solid thin line) corresponds to a not-so-pushy $\msem=0.23137\cdot 10^{-2}$  ($r=20$) that, nevertheless, is already quite effective to identify and remove the outliers.
The performance in executing the tasks is corroborated by the right diagram which illustrates the history of the maximal-entropy weights when the reduction factor increases.
The outliers' weights $p_{5}$ and $p_{10}$ decrease rapidly; in particular, when $r=20$, they become $p_5\approx 8.4\cdot 10^{-9}$ and $p_{10}\approx 2.2\cdot 10^{-3}$.
The other weights remain close to the initial value $p_i(r=1)=1/n=1/20=0.05$;  
we have included in the diagram only the curves for {\blue $p_{1}$ and $p_{20}$ as} reference because the remaining weights' curves lie in between. 
Yet the weights of the perturbed data ($i\not=5,10$) do not assume exactly the constant value $1/18=0.0556$ and, therefore, the mew-parabola should be formally distinct from the dashed uw-parabola obtained after manual removal of the outliers from the data set.
Nonetheless, the two objects are essentially equivalent as evidenced by the left diagram.
In conclusion, MEM identifies and removes outliers automatically, without the need of any sort of \textit{human} intervention.

\subsection{Monna Lisa's eyes}
The conclusion of \Rse{subsec:out-rem} suggests another interesting application of MEM in the field of image processing.
As final test case, we consider the problem of removing noise, such as that produced by dust artifacts caused by contaminants on digital-camera lenses or by dead, stuck or defective pixels in sensors, from a digital image.  
To this end, we introduce an artificial Gaussian distributed noise into an original grayscale picture and test the performance of MEM to recover the original image. 
The original grayscale image is stored in a $n \times m$ matrix $\Pi$ whose element $\pi_{ij}$ carries the  $(i,j)$-pixel's intensity between 0 for black and 1 for white.
The experiment proceeds through the following key steps:
\begin{enumerate}
   \item  {\it Selection of pixels that will be affected by noise} \\
          We generate a $n \times m$ random logical matrix $\Omega$ with elements 
          \begin{equation}
             \omega_{ij} = 
          	   \begin{cases} \; 1 & \quad \text{with probability $\mbox{\textsf{P}}$}  \\[1.5ex] 
          	                 \; 0 & \quad \text{with probability $1-\mbox{\textsf{P}}$}  \\[0ex] 
          	   \end{cases} \label{omega}
          \end{equation}
          for a given $\mbox{\textsf{P}}\in[0,1]$.
   \item  {\it Introduction of noise in the original image.} \\
          We generate the perturbed matrix $\widetilde\Pi = \Pi + \delta\,\Pi$ by adding a (truncated) Gaussian noise only to those pixels filtered by the mask $\Omega$. 
          Thus, we set: 
		  \begin{equation}\label{pert}
		     \subeqn{\tilde{\pi}_{ij} = 
		  	   \begin{cases} \; \min\{1,\max\{0,\pi_{ij}+\alpha \eta_{ij}\}\}  & \quad \text{if $\omega_{ij}=1$}  \\[1.5ex] 
		  	                 \; \pi_{ij}                                       & \quad \text{if $\omega_{ij}=0$}  \\[0ex] 
		  	   \end{cases}}{i=1,\ldots,n\,;\,j=1,\dots,m}{} 
		  \end{equation}
		  In \Req{pert}, \mbox{$\eta_{ij} \sim {\cal N}(0,1)$} is a random variable  distributed normally with mean $0$ and variance $1$, $\alpha\in[0,1]$ is a safety factor meant to increase the probability that $\pi_{ij}+\alpha \eta_{ij}$ belongs to the interval $[0,1]$; if this is not the case the $\min$ and $\max$ functions come into play by forcing $\tilde{\pi}_{ij}$ to be $0$ or $1$ according to whether $\pi_{ij}+\alpha \eta_{ij}$ is negative or greater than $1$, respectively. 
   \item  {\it Iteration procedure to remove noise.} \\
          We apply our FM algorithm iteratively to rows and/or columns of the matrix $\widetilde\Pi$ by progressively reducing the value of $\msem$ and intervening on those pixels whose associated weights become negligible by replacing their intensity values with those predicted by the maximal-entropy approximation. 
\end{enumerate}
Clearly, the aim of steps 1 and 2 is to prepare the noisy image. 
We ought to emphasize here that, due to the random nature of the noise intensity, the detection of noisy pixels' coordinates is much more involved than it would be when dealing with situations in which noisy pixels have an intensity known a priori; as example of the latter situation, we mention the \textit{salt and pepper} noise that consists only of either the minimum ($0$) or maximum ($1$) intensity values. 
In our case study, the intensity of corrupted pixels may assume any value in the interval $[0,1]$ and when that intensity is close to the intensity of neighboring pixels, the identification of noisy pixels is obviously a difficult task to handle.  
The detection of corrupted pixels and the de-noising technique are carried out concurrently and incorporated in a dynamic procedure. 
For a given prescribed value of $\msem$, a noisy pixel will be perceived as an outlier and correctly de-noised. 
In other words, $\msem$ acts as a trimmer to adjust the sensitivity of the method in deciding whether or not a given pixel should be labeled as corrupted. 
It may occur that noise-free pixels could be erroneously labeled as corrupted if $\msem$ is chosen too small; in order to prevent false positives, we have iteratively reduced  $\msem$ until a satisfactory balance between sensitivity and specificity was attained.  
\begin{figure}[h]
   \centering
   \resizebox{.325\textwidth}{!}{\includegraphics*[trim=0 0 0 0]{\figdir/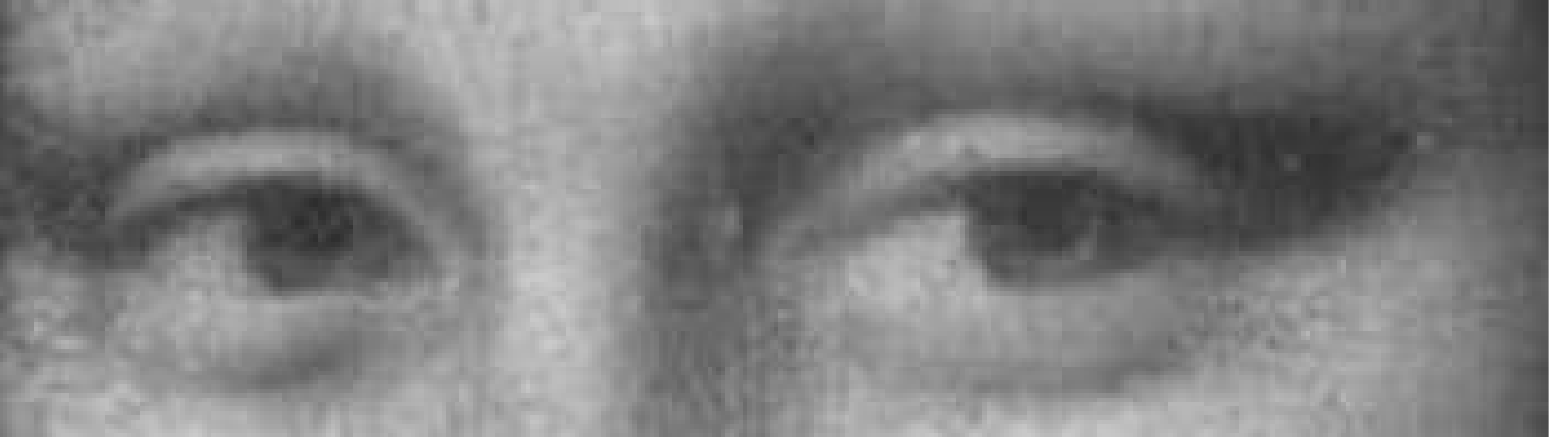}}\hfill      
   \resizebox{.325\textwidth}{!}{\includegraphics*[trim=0 0 0 0]{\figdir/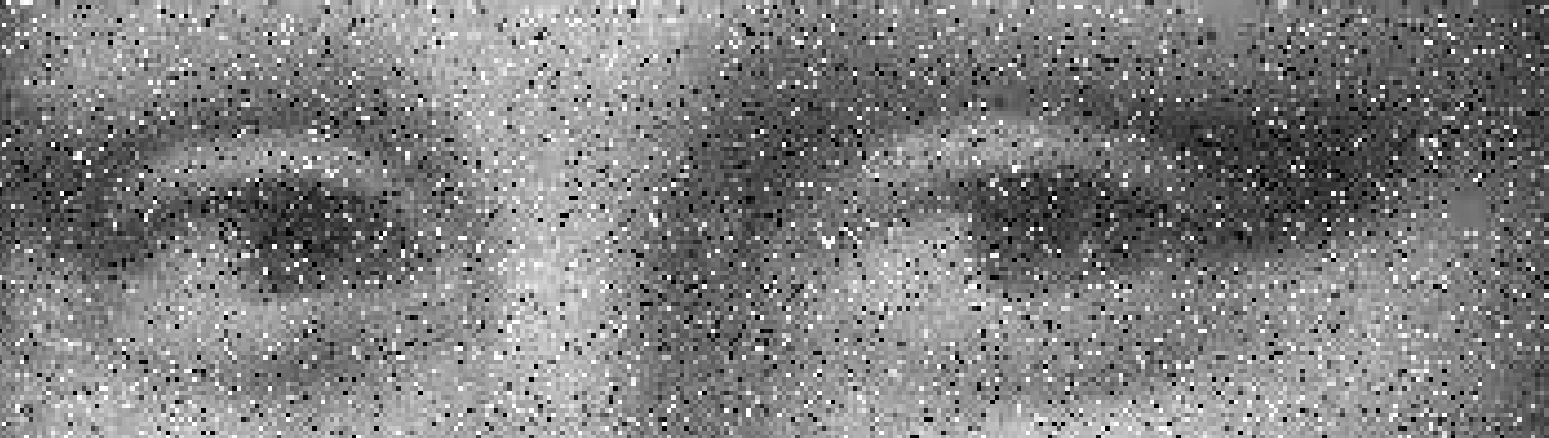}}\hfill      
   \resizebox{.325\textwidth}{!}{\includegraphics*[trim=0 0 0 0]{\figdir/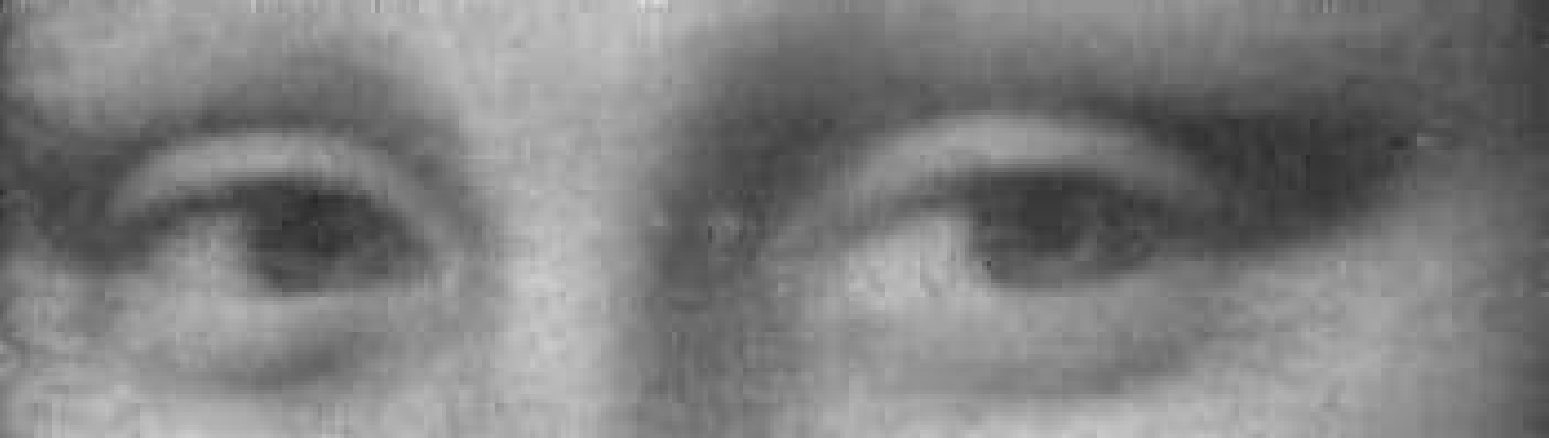}}\hfill      
   \caption{The three stages of the image-processing experiment on a close-up of Monna Lisa's eyes. 
            Left:   original image. 
            Center: after random-noise introduction. 
            Right:  after cleaning by iterative application of the FM algorithm.\label{pic-monnalisa}}   
\end{figure}

We have applied the procedure described above in steps 1--3  to the image shown in \Rfi{pic-monnalisa} (left), a grayscale close-up of Monna Lisa's eyes, stored in a matrix $\Pi$ with dimension $n \times m = 99 \times 350$. 
In step 1, we have set $\mbox{\textsf{P}}=0.15$ which implies a noise density of about 15\%; in step 2, we have set $\alpha=1/2$. 
The resulting noisy image is shown in \Rfi{pic-monnalisa} (center). 
Then we have executed the FM algorithm iteratively by sweeping columns and rows of the matrix $\widetilde\Pi$. 
Starting from the initial value corresponding to uniform weights, $\msem$ was progressively reduced and, at each iteration, the pixel $(i,j)$ was labeled {\it corrupted} if its weight dropped below a given tolerance. 
The resulting cleaned image is shown in \Rfi{pic-monnalisa} (right) and looks satisfactorily good (we dare to say). 

\section{Conclusions}

The experience we have acquired with the test cases described in this work, and additional others that we have not included, gives us strong confidence in the good performance of MEM.
Its application is relatively straightforward and we are amazed of how much we have been able to deal with just with the use of polynomials as approximating functions and one-dimensional data sets.
In fact, we are absolutely convinced that we have just scratched the surface; the full potentiality of the method lies still undiscovered and awaiting to be brought to fruition.
In this regard, future work will definitely deal with the necessity to extend the MEM we have presented in here to include within its reach {\blue interesting and useful features such as} the use of approximating functions other than polynomials, the presence of uncertainties for both $Y_{i}$ and $X_{i}$ data, {\blue handling} multidimensional data sets, and, last but not least, {\refereetwo to investigate the effective implementation of the method in terms of computational complexity and stability of the numerical algorithm.} 
Moreover, we consider important to study the mathematical properties of the entropy, its first and second derivatives as functions of the minimal mean squared error; regarding this last aspect, we feel confident that much analogy can de drawn from statistical thermodynamics {\blue if one yields at the scientific-curiosity push to investigate matters like, for example, whether or not the thermodynamic-stability body of knowledge can be exported to and reinterpreted in the MEM mathematical context on the basis of the analog connection between the role of energy in statistical thermodynamics and the role of mean squared error in MEM. A fascinating problem indeed!}



\bibliography{ms} 

\end{document}